%% file: authsamp.tex
\documentclass[graybox]{svmult}

\usepackage{type1cm}        % activate if the above 3 fonts are
\usepackage{makeidx}         % allows index generation
\usepackage{graphicx}        % standard LaTeX graphics tool
\usepackage{multicol}        % used for the two-column index
\usepackage[bottom]{footmisc}% places footnotes at page bottom

\usepackage{newtxtext}       % 
\usepackage{newtxmath}       % selects Times Roman as basic font

\makeindex             % used for the subject index
\begin{document}

\title*{BDDC preconditioners for a space-time finite element discretization
  of parabolic problems}  
\titlerunning{BDDC for STFEM of parabolic problems} %for an abbreviated version of
% your contribution title if the original one is too long
\author{Ulrich Langer and Huidong Yang}
% Use \authorrunning{Short Title} for an abbreviated version of
% your contribution title if the original one is too long
\institute{Ulrich Langer \at Johann Radon Institute, Altenberg Strasse 69, 
  4040 Linz, Austria, \email{ulrich.langer@ricam.oeaw.ac.at} 
\and Huidong Yang \at Johann Radon Institute, Altenberg Strasse 69, 
  4040 Linz, Austria, \email{huidong.yang@oeaw.ac.at}}
%
% Use the package "url.sty" to avoid
% problems with special characters
% used in your e-mail or web address
%
\maketitle

\abstract{% This paper deals with balanced domain decomposition by constraints
  %(BDDC) method for a space-time finite element discretization of parabolic
  %problems, where time is considered as another spatial coordinate and the
  %finite elements are continuous and piecewise linear simultaneously in space
  %and time. We consider BDDC preconditioned GMRES methods for solving the discrete
  %space-time finite element Schur complement equations on the
  %interface. Numerical studies demonstrate robustness of the preconditioners
  %to some extent.
  This paper deals with balanced domain decomposition by constraints
  (BDDC) method for solving 
  %the 
  large-scale linear systems of algebraic equations
  arising from the space-time finite element discretization of parabolic
  initial-boundary value problems.
  The time is considered as just another spatial coordinate, and the
  finite elements are
  %chosen to be 
  continuous and piecewise linear
  on unstructured simplicial space-time meshes. 
  We consider BDDC preconditioned GMRES methods for solving the 
  %discrete
  space-time finite element Schur complement equations on the
  interface. Numerical studies demonstrate robustness of the preconditioners
  to some extent.}

\section{Introduction}\label{sec:intro}
Continuous space-time finite element methods for parabolic problems have been
recently studied, e.g.,  in \cite{REB17,UL16,UL18,OS15}. 
%As 
The main common features of these methods
%, that are distinguished from conventional time stepping methods, 
are very different from those of time-stepping methods.
Time is considered to be just another spatial coordinate. 
The variational formulations are studied in the full space-time cylinder 
that is then decomposed into arbitrary admissible simplex elements. 
In this work, we follow the space-time finite element discretization scheme
proposed in 
\cite{UL18} 
%\cite{UL18,AS17} 
for a model initial-boundary value problem, using continuous
and piecewise linear finite elements  in space and time simultaneously. 

%On the other hand, it 
It is a challenge task to efficiently solve the large-scale linear
system of algebraic equations arising from the space-time finite element
discretization of parabolic problems. In this work, as a preliminary study,
we use the balanced domain decomposition by constraints (BDDC \cite{CD03,JM03,JM05}) 
preconditioned GMRES method to solve this system efficiently.
%, which aims to save the computational time. 
We mention that robust preconditioning for space-time
isogeometric analysis schemes for parabolic evolution problems 
has been reported in \cite{CH1802,CH18}.

The remainder of the paper is organized as follows: Sect.~\ref{sec:stfem}
deals with the space-time finite element discretization for 
%a model parabolic
a parabolic model 
problem. In Sect.~\ref{sec:stbddc}, we discuss BDDC preconditioners that are
used to solve the 
%discrete 
linear system of algebraic equations. 
Numerical results are shown and discussed in Sect.~\ref{sec:numa}. 
Finally, some conclusions are drawn in Sect.~\ref{sec:cons}.

\section{The space-time finite element discretization}
\label{sec:stfem}
%\subsection{A parabolic model equation}\label{subsec:model}
The following parabolic initial-boundary value problem is considered as our model problem: Find
$u:\overline{Q}\rightarrow{\mathbb R}$ such that
\begin{equation}\label{eq:modelpar}
\partial_t u - \Delta_x u  =  f
\textup{ in } Q,\quad
u=0 \textup{ on }\Sigma,\; 
u=u_0\textup{ on }\Sigma_0,
\end{equation}
where $Q := \Omega \times (0,T)$, 
$\Omega\subset{\mathbb R}^2$ is a sufficiently smooth and bounded computational domain, 
%$d=2$, 
$\Sigma := \partial \Omega \times (0,T)$, $\Sigma_0:=\Omega\times\{0\}$,
$\Sigma_T:=\Omega\times\{T\}$. 

%\subsection{A space-time variational formulation}\label{subsec:stvar}
Let us now introduce the following Sobolev spaces: 
\begin{equation*}
  \begin{aligned}
    %H^{l,k}(Q)&=\{u\in L_2(Q) : \partial_x^\alpha u\in L_2(Q), 0\leq|\alpha|\leq l,\partial_t^i u \in L_2(Q), i=0,...,k\},\\
    H_0^{1,0}(Q)&=\{u\in L_2(Q):\nabla_xu\in[L_2(Q)]^2, u=0\textup{ on }
    \Sigma\},\\
    H_{0,\bar{0}}^{1,1}(Q)&=\{u\in L_2(Q):
    \nabla_xu\in[L_2(Q)]^2, \partial_tu\in L_2(Q) \textup{ and } u|_{\Sigma\cup\Sigma_T}=0\}, \\
    H_{0,\underline{0}}^{1,1}(Q)&=\{u\in L_2(Q):
    \nabla_xu\in[L_2(Q)]^2, \partial_tu\in L_2(Q) \textup{ and } u|_{\Sigma\cup\Sigma_0}=0\}. \\
  \end{aligned}
\end{equation*}
Using the classical approach \cite{OAL85,OAL68}, the variational formulation
for the 
%model parabolic 
parabolic model 
problem (\ref{eq:modelpar}) reads as follows: Find
$u\in H_0^{1, 0}(Q)$ such that  
\begin{equation}\label{eq:var}
a ( u , v ) = l ( v ), \quad \forall v\in H_{0,\bar{0}}^{1,1}(Q),
\end{equation}
where 
\begin{equation*}
  \begin{aligned}
    a(u,v)&=-\int_Q u(x,t)\partial_tv(x,t)dxdt+
    \int_Q\nabla_x u(x,t)\cdot\nabla_x u(x,t)dxdt,\\
    l(v)&=\int_Q f(x,t)v(x,t)dxdt + \int_\Omega u_0(x)v(x,0)dx.
  \end{aligned}
\end{equation*}
\begin{remark}[Parabolic solvability and regularity \cite{OAL85,OAL68}]
  If $f\in L_{2,1}(Q):=\{v: \int_0^T\|v(\cdot, t)\|_{L_2(\Omega)}dt<\infty\}$
  and $u_0\in L_2(\Omega)$, then there exists a unique generalized solution 
  $u\in H_0^{1,0}(Q)\cap V_{2}^{1,0}(Q)$ of 
  %the variational problem,
  (\ref{eq:var}),
  %{\color{red} where $V_{2}^{1,0}(Q) = ... $.}
  where $V_{2}^{1,0}(Q):=\{u\in H^{1,0}(Q):|u|_Q<\infty\textup{ and }
  \displaystyle\lim_{\Delta t\rightarrow 0}\|u(\cdot , t+\Delta t) -
  u(\cdot , t)\|_{L_2(\Omega)}=0,\textup{ uniformly on }[0, T]\}$,
  and 
    $|u|_Q := \displaystyle\max_{0\leq\tau\leq
    t}\|u(\cdot,\tau)\|_{L_2(\Omega)}+\|\nabla_x u\|_{L_2(\Omega\times(0,t))}$. 
  If $f\in L_2(Q)$ and $u_0\in H_0^{1}(\Omega)$, then the generalized solution 
  $u$ belongs to $H_0^{\Delta,1}(Q):=\{v\in H_0^{1,1}(Q) : \Delta_x u\in L_2(Q) \}$ and
  continuously depends on $t$ in the norm of the space $H_0^1(\Omega)$.
  %, where $\Delta_xu\in L_2(Q)$ means $\nabla_x u\in H(\textup{div}; Q)$.
  \label{rm:reg} 
\end{remark}
%\subsection{Discrete space-time variational formulation}\label{subsec:disstfem}
To derive the 
%discrete variational formulation, 
space-time finite element scheme,
we mainly follow the 
%scheme
approach 
proposed in \cite{UL18}.
%\cite{UL18,AS17}. 
Let 
$V_h=\textup{span}\{\varphi_i\}$ be the span of continuous and piecewise
linear basis functions 
$\varphi_i$ on shape regular finite elements of an admissible triangulation
${\mathcal T}_h$. 
Then we define $V_{0h}=V_h\cap H_{0,\underline{0}}^{1,1}(Q)=\{v_h\in V_h:
v_h|_{\Sigma\cup\Sigma_0}=0\}$. For convenience, we consider homogeneous
initial conditions, i.e., $u_0=0$ on $\Omega$. Multiplying the PDE $\partial_t
u  - \Delta_xu = f$ on $K\in{\mathcal T}_h$ by an element-wise time-upwind
test function $v_h+\theta_K h_K\partial_t v_h$, $v_h\in V_{0h}$, we get 
\begin{equation*}
  \begin{aligned}
    &\int_K(\partial_t u v_h + \theta_K h_K \partial_t u \partial_t v_h
    -\Delta_x u ( v_h + \theta_K h_K \partial_t v_h))dxdt = \\
    &\int_K f (v_h + \theta_K h_K \partial v_h)dxdt.
  \end{aligned}
\end{equation*}
Integration by parts (the first part) with respect to the space and summation lead to 
\begin{equation*}
  \begin{aligned}
    &\displaystyle\sum_{K\in{\mathcal T}_h}\int_K (\partial_t u v_h 
    + \theta_K h_K \partial_t u \partial_tv_h + 
    \nabla_x u \cdot \nabla_x v_h - \theta_K h_K \Delta_x u \partial_t v_h )dxdt\\
    &-\displaystyle\sum_{K\in{\mathcal T}_h}\int_{\partial K} n_x\cdot
    \nabla_x u v_h ds 
    = \displaystyle\sum_{K\in{\mathcal T}_h}\int_K f (v_h + \theta_K h_K \partial_t v_h) dxdt.
  \end{aligned}
\end{equation*}
Since $n_x\cdot\nabla_x u$ is continuous on inner boundary of $K$, $n_x=0$ on
$\Sigma_0\cup\Sigma_T$, and $v_h=0$ on $\Sigma$, the term
$-\sum_{K\in{\mathcal T}_h}\int_{\partial K} n_x\cdot\nabla_x u
v_h ds$ vanishes. 

%Let $H_{0,\underline{0}}^{2,1}({\mathcal T}_h):=\{v\in
%H_{0,\underline{0}}^{1,1}(Q) : v|_K \in H^{2,1}(K), \forall K\in{\mathcal T}_h\}$. 
%For the solution $u\in H_{0,\underline{0}}^{2,1}({\mathcal T}_h)\cap
%H_0^{\Delta, 1}(Q)$ of the Ladyzhenskaya's variational formulation, 
%we have the consistency identity, 
%
%Let $H_{0,\underline{0}}^{\Delta,1}({\mathcal T}_h)
%:=\{v\in H_{0,\underline{0}}^{1,1}(Q) : \Delta_x v|_K \in L_2}(K), \forall K\in{\mathcal T}_h\}$. 
If the solution $u$ of 
%the variational formulation 
%{\color{red} (LABEL) } 
{(\ref{eq:var})} 
belongs to 
$H_{0,\underline{0}}^{\Delta,1}({\mathcal T}_h):= \{v\in
H_{0,\underline{0}}^{1,1}(Q) : \Delta_x v|_K \in L_2(K), \forall K\in{\mathcal
  T}_h\}$, 
%cf. Remark~{\color{red} (LABEL) }, 
cf. Remark~\ref{rm:reg}, 
then the consistency identity
\begin{align}
  a_h(u, v_h) = l_h (v_h),\quad v_h\in V_{0h},
\end{align}
holds, where 
\begin{align*}
  a_h(u, v_h)&:=\displaystyle\sum_{K\in{\mathcal T}_h}\int_K (\partial_t u v_h + \theta_K h_K \partial_t
  u\partial_t v_h  + \nabla_x u\cdot \nabla_x v_h - \theta_K h_K \Delta_x u \partial_t v_h)dxdt,\\
  l_h(v_h)&:=\displaystyle\sum_{K\in{\mathcal T}_h}\int_K f(v_h+\theta_K h_K\partial_t v_h)dxdt.
\end{align*}
With the restriction of the solution to the 
finite-dimensional
subspace $V_{0h}$, 
the 
%discrete space-time variational formulation 
space-time finite element scheme reads as follows: Find $u_h\in V_{0h}$ such that 
\begin{align}
  a_h(u_h, v_h) = l_h (v_h),\quad v_h\in V_{0h}.
\end{align}
%where 
%\begin{align*}
%  a_h(u_h, v_h)&=\displaystyle\sum_{K\in{\mathcal T}_h}\int_K (\partial_t u_h v_h + \theta_K h_K \partial_t
%  u_h\partial_t v_h  + \nabla_x u_h\cdot \nabla_x v_h - \theta_K h_K\Delta_x u_h \partial_t v_h )dxdt,\\
%  l_h(v_h)&=\displaystyle\sum_{K\in{\mathcal T}_h}\int_K f(v_h+\theta_K h_K\partial_t v_h)dxdt.
%\end{align*}
%It is easy to see that 
Thus,
we have the Galerkin orthogonality: $a_h(u -u_h, v_h)=0$, 
$\forall v_h\in V_{0h}$. 
\begin{remark}
  Since we use continuous and piecewise linear trial functions, the integrand 
  $-\theta_K h_K\Delta_x u_h \partial_t v_h$ vanishes element-wise, which simplifies the
  implementation. 
\end{remark}
\begin{remark}
  On fully unstructured meshes, $\theta_k=O(h_k)$ \cite{UL18}; on uniform
  meshes, $\theta_k=\theta=O(1)$ \cite{UL16}. In this work, we have used
  $\theta=0.5$ and $\theta=2.5$ on uniform meshes for testing robustness of
  the BDDC preconditioners. 
\end{remark}

%\subsection{The a priori discretization error estimation}
It was shown in \cite{UL18} that the 
%discrete 
bilinear form 
$a_h(\cdot,\cdot)$
is
$V_{0h}$-coercive: $a_h (v_h, v_h)\geq \mu_c\|v_h\|_h^2$, $\forall v_h\in V_{0h}$ with
respect to the norm $\|v_h\|_h^2=\sum_{K\in{\mathcal
    T}_h}(\|\nabla_xv_h\|_{L_2(K)}^2 + \theta_K h_K \|\partial_t
v_h\|_{L_2(K)}^2)+ \frac{1}{2} \|v_h\|^2_{L_2(\Sigma_T)}$. 
Furthermore, the bilinear form is bounded on $V_{0h,*}\times V_{0h}$:
$|a_h(u,v_h)|\leq \mu_b\|u\|_{0h,*}\|v_h\|_h$, $\forall u\in V_{0h,*}$,
$\forall v_h\in V_{0h}$, where 
%$V_{0h,*}=H_{0}^{\Delta, 1}(Q)\cap H^2 ({\mathcal T}_h)+V_{0h}$ 
$V_{0h,*}=H_{0,\underline{0}}^{\Delta, 1}({\mathcal T}_h)+V_{0h}$
equipped with the norm 
%$\|v\|_{0h,*}^2=\|v\|_h^2 + \sum_{K\in{\mathcal T}_h}(\theta_K h_K)^{-1}\|v\|_{L_2(K)}^2$ $+
%\sum_{K\in{\mathcal T}_h}\theta_k h_k |v|_{H^2(K)}^2$.
$\|v\|_{0h,*}^2=\|v\|_h^2 + \sum_{K\in{\mathcal T}_h}(\theta_K
h_K)^{-1}\|v\|_{L_2(K)}^2$ $+
\sum_{K\in{\mathcal T}_h}\theta_k h_k \|\Delta_x v\|_{L_2(K)}^2$.
%
%We now define the broken Sobolev space 
%$H^2({\mathcal T}_h):=\{ v\in L_2(Q) : v|_K \in H^2(K)$ $\forall K\in{\mathcal T}_h\}$ 
%equipped with the broken Sobolev semi-norm $|v|_{H^2 ({\mathcal T}_h)}:=\sum_{K\in {\mathcal T}_h}|v|_{H^2(K)}^2$.
%Using a proper interpolation operator $\Pi_h$ mapping 
%from $H_{0,\underline{0}}^{1,1}\cap H^k(Q)$ to $V_{0h}$ for $k>\frac{d+1}{2}$, 
%then $\|u-u_h\|_h \leq \|u- \Pi_h u\|_h + \|\Pi_h u -u_h\|_h$. 
%The first term $\|u- \Pi_h u\|_h$ can be bounded by using the interpolation
%error estimate, and the second term $\|\Pi_h u -u_h\|_h$ by using ellipticity,
%Galerkin orthogonality and boundedness of the bilinear form. The discretization
%error estimate is (see more details in \cite{UL18}):
%$\|u-u_h\|_h \leq C (\sum_{K\in{\mathcal T}_h}h_K^2|u|^2_{ H^2 (K)})^{1/2}$, 
%for the solution 
%$u\in H_{0,\underline{0}}^{1,1}\cap H^k(Q) \cap H^2({\mathcalT}_h)$ 
%$and $u_h\in V_{0h}$, with $C>0$, independent of mesh size.
%
%
%{\color{blue}
Let $l$ and $k$ be positive 
%integers 
reals
such that $l\geq k> 3/2$.
%$s=\min\{l,2\}$. UL: ???
We now define the broken Sobolev space $H^s({\mathcal
  T}_h):=\{ v\in L_2(Q) : v|_K \in H^s(K)$ $\forall K\in{\mathcal T}_h\}$
equipped with the broken Sobolev semi-norm 
$|v|^2_{H^s ({\mathcal T}_h)}:=\sum_{K\in {\mathcal T}_h}|v|_{H^s(K)}^2$. 
Using the Lagrangian interpolation operator $\Pi_h$ mapping 
%from 
$H_{0,\underline{0}}^{1,1}(Q)\cap H^k(Q)$ to $V_{0h}$, 
%then 
we obtain 
$\|u-u_h\|_h \leq \|u- \Pi_h u\|_h + \|\Pi_h u -u_h\|_h$. 
The 
%first 
term $\|u- \Pi_h u\|_h$ can be bounded by means of the interpolation
error estimate, and the 
%second 
term $\|\Pi_h u -u_h\|_h$ by using ellipticity,
Galerkin orthogonality and boundedness of the bilinear form. The discretization
error estimate
%\cite{UL18}, 
$\|u-u_h\|_h \leq C (\sum_{K\in{\mathcal T}_h}h_K^{2(l-1)}|u|^{2}_{ H^{l} (K)})^{1/2}$
holds for the solution $u$ 
%$u\in H_{0,\underline{0}}^{1,1}(Q)\cap H^k(Q) \cap H^s({\mathcal T}_h)$ 
provided that $u$ belongs to $H_{0,\underline{0}}^{1,1}(Q)\cap H^k(Q) \cap H^l({\mathcal T}_h)$,
and the finite element solution $u_h\in V_{0h}$, 
%with a positive constant $C>0$, independent of mesh size; see \cite{UL18}. 
where $C>0$, independent of mesh size; see \cite{UL18}. 
%$C$ is a  positive constant that is independent of mesh size; see \cite{UL18}. 
%The assumption $k>3/2$ can be further relaxed; see more details in \cite{UL18}.
%} 

\section{Two-level BDDC preconditioners}\label{sec:stbddc}
%
%The linear system of algebraic equations 
After the space-time finite element
discretization of the model problem (\ref{eq:modelpar}),
the linear system of algebraic equations 
reads as follows:  
%\begin{align}
%  &Kx=f,\\
%  &K:=
%  \begin{bmatrix}
%    K_{II} & K_{I\Gamma} \\
%    K_{\Gamma I} & K_{\Gamma\Gamma} \\
%  \end{bmatrix}, \quad
%  x:=
%  \begin{bmatrix}
%    x_I\\
%    x_\Gamma \\
%  \end{bmatrix},\quad
%  f:=
%  \begin{bmatrix}
%    f_I\\
%    f_\Gamma\\
%  \end{bmatrix}, \quad
%  K_{II}=\text{diag}\left[K_{II}^1,...,K_{II}^{N}\right], \nonumber
%\end{align}
\begin{align}
\label{eq:Ku=f}
  Kx=f,
\end{align}
with $K:=
\begin{bmatrix}
  K_{II} & K_{I\Gamma} \\
  K_{\Gamma I} & K_{\Gamma\Gamma} \\
\end{bmatrix}$, 
$x:=
\begin{bmatrix}
  x_I\\
  x_\Gamma \\
\end{bmatrix}$, 
$f:=\begin{bmatrix}
  f_I\\
  f_\Gamma\\
\end{bmatrix}$, 
$K_{II}=\text{diag}\left[K_{II}^1,...,K_{II}^{N}\right]$, 
where $N$ denotes the number of polyhedral subdomains
$\Omega_i$ from a non-overlapping domain decomposition of $\Omega$. 
In system (\ref{eq:Ku=f}), we have decomposed the degrees of freedom into the ones associated with
the internal ($I$) and interface ($\Gamma$) nodes, respectively. We aim to
solve the 
%reduced 
Schur-complement
system living on the interface:
%\begin{align}\label{eq:schur}
%  Sx_\Gamma &=g_\Gamma, \\
%    S&:=K_{\Gamma\Gamma} - K_{\Gamma I}K_{II}^{-1}K_{I\Gamma},\quad
%    g:=f_\Gamma - K_{\Gamma I}K_{II}^{-1}f_I. \nonumber
%\end{align}
\begin{align}\label{eq:schur}
  Sx_\Gamma &=g_\Gamma,
\end{align}
with $S:=K_{\Gamma\Gamma} - K_{\Gamma I}K_{II}^{-1}K_{I\Gamma}$ and 
$g:=f_\Gamma - K_{\Gamma I}K_{II}^{-1}f_I$. 

Following \cite{JM05} (see also details in \cite{LJ06}), Dohrmann's
(two-level) BDDC preconditioners $P_{BDDC}$ for the interface Schur
complement equation (\ref{eq:schur}), originally proposed 
for symmetric and positive definite systems in \cite{CD03,JM03}, 
%may 
can 
be written in the form 
\begin{equation}
  P_{BDDC}^{-1}=R_{D,\Gamma}^T(T_{sub}+T_0)R_{D,\Gamma},
\end{equation}
where the scaled operator $R_{D,\Gamma}$ is the direct sum of restriction
operators $R_{D,\Gamma}^i$ mapping the global interface vector to its component on
local interface $\Gamma_i:=\partial\Omega_i\cap \Gamma$, with a proper scaling
factor. 

Here the coarse level correction operator $T_0$ is constructed  as 
\begin{equation}
  T_0 = \Phi (\Phi^T S \Phi)^{-1} \Phi^T
\end{equation}
with the coarse level basis function matrix 
$\Phi=\left[(\Phi^1)^T,\cdots, (\Phi^N)^T\right]^T$, where the basis function
matrix $\Phi^i$ on each subdomain interface is obtained by solving the
following augmented system:
\begin{equation}
  \begin{aligned}
    \begin{bmatrix}
      S^i & \left(C^i\right)^T \\
      C^i & 0
    \end{bmatrix}
    \begin{bmatrix}
      \Phi^i\\
      \Lambda^i\\
    \end{bmatrix}
    =
    \begin{bmatrix}
      0\\
      R_{\Pi}^i\\
    \end{bmatrix}.
  \end{aligned}
\end{equation}
with the given primal constraints $C^i$ of the subdomain $\Omega_i$ and
the vector of Lagrange multipliers on each column of $\Lambda^i$. The number
of columns of each $\Phi^i$ equals to the number of global coarse level
degrees of freedom, typically living on the subdomain corners, and/or
interface edges, and/or faces. 
Here the restriction operator $R_{\Pi}^i$ maps
the global interface vector in the continuous primal variable space on the
coarse level to its component on $\Gamma_i$.  

The subdomain correction operator $T_{sub}$ is defined as 
\begin{equation}
  T_{sub}=\displaystyle\sum_{i=1}^N\left[(R_\Gamma^i)^T\quad 0\right]
  \begin{bmatrix}
    & S^i & \left(C^i\right)^T\\
    & C^i & 0
  \end{bmatrix}^{-1}
  \begin{bmatrix}
    & R_\Gamma^i\\
    &0
  \end{bmatrix},
\end{equation}
with vanishing primal variables on all the coarse levels. Here the restriction
operator $R_\Gamma^i$ maps global interface vectors to their components on $\Gamma_i$.

\section{Numerical experiments}\label{sec:numa}
\begin{figure}[t]
\sidecaption
\includegraphics[scale=.15]{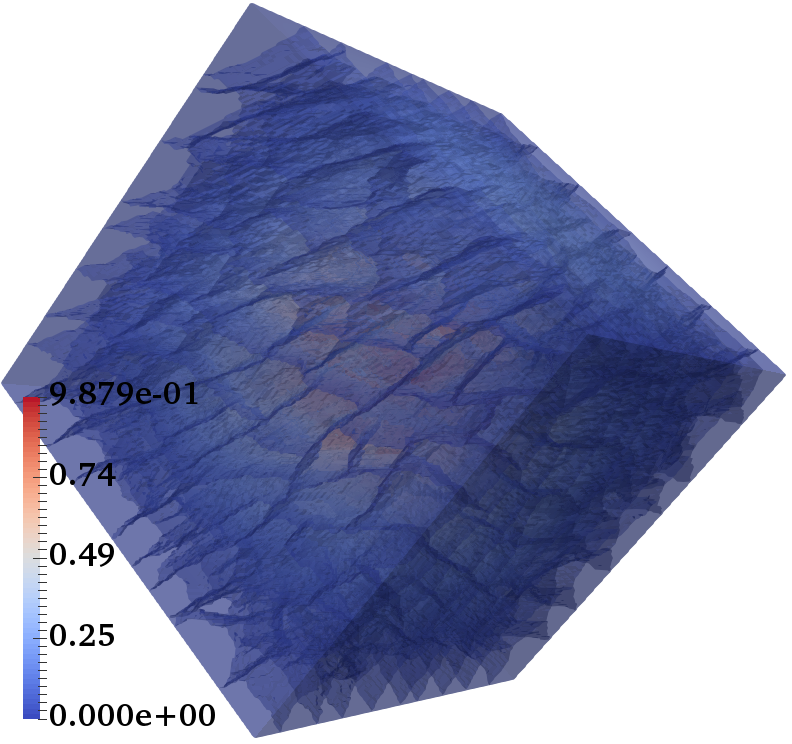}\quad
\includegraphics[scale=.15]{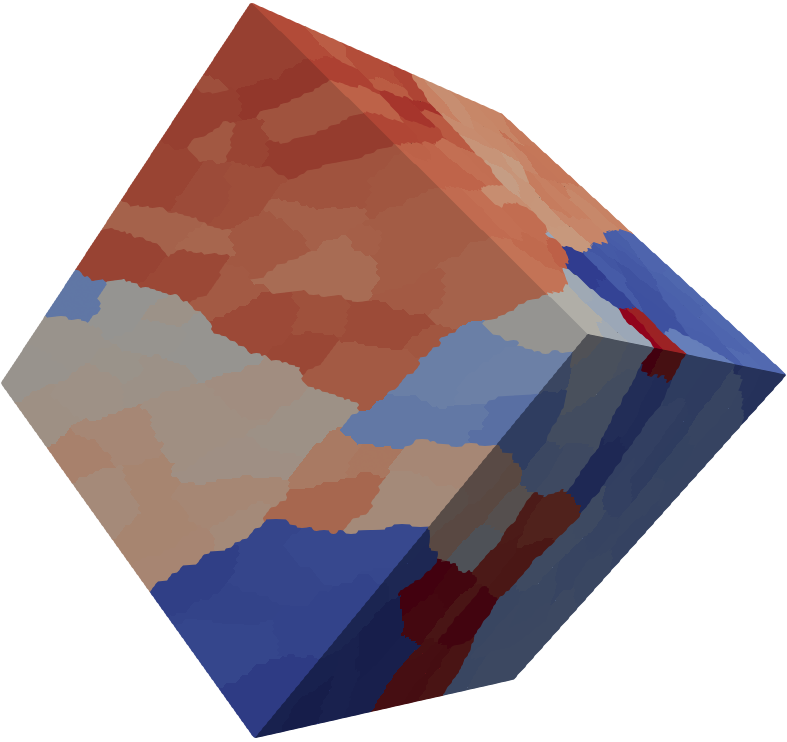}
\caption{Solution (left), space-time domain decomposition (right) with
  $129^3$ degrees of freedom and $512$ subdomains.} 
\label{fig:soldd}
\end{figure}

%\subsection{Convergence}
%We consider an exact solution $u(x,y,t)=\sin(\pi x)\sin(\pi y)\sin(\pi t)$, $Q=(0,1)^3$,
%\begin{table}
%  \begin{center}
%    \begin{tabular}{ crcccc }
%      \hline
%      $\theta=0.5$ & \#Dofs & $\|u-u_h\|_{L_2(0,T; H_0^1(\Omega))}$ & eoc  \\ \hline
%      &$4,913$ $(17\times 17\times 17)$ & $3.7404\cdot 10^{-1}$ & $-$ \\ \hline
%      &$35,937$ $(33\times 33\times 33)$ & $1.9769\cdot 10^{-1}$ & $0.9200$
%      \\ \hline 
%      &$274,625$ $(65\times 65\times 65)$ & $1.0168\cdot 10^{-1}$ & $0.9592$ \\ \hline
%      &$2,146,689$ $(129\times 129\times 129)$ & $5.1600\cdot 10^{-2}$ & $0.9786$ \\ \hline
%    \end{tabular}
%  \end{center}
%\end{table}
%\begin{table}
%  \begin{center}
%    \begin{tabular}{ crcccc }
%      \hline
%      $\theta=2.5$ & \#Dofs & $\|u-u_h\|_{L_2(0,T; H_0^1(\Omega))}$ & eoc  \\ \hline
%      &$4,913$ $(17\times 17\times 17)$ & $1.2093\cdot 10^{-0}$ & $-$ \\ \hline
%      &$35,937$ $(33\times 33\times 33)$ & $6.9825\cdot 10^{-1}$ & $0.7924$
%      \\ \hline
%      &$274,625$ $(65\times 65\times 65)$ & $3.8213\cdot 10^{-1}$ & $0.8697$ \\ \hline
%      &$2,146,689$ $(129\times 129\times 129)$ & $2.0183\cdot 10^{-1}$
%      & $0.9209$ \\ \hline
%    \end{tabular}
%  \end{center}
%\end{table}
We use 
%an exact solution 
$u(x,y,t)=\sin(\pi x)\sin(\pi y)\sin(\pi t)$ as exact solution of (\ref{eq:modelpar})
in $Q=(0,1)^3$; see the left
plot in Fig.~\ref{fig:soldd}. We
perform uniform mesh refinements of $Q$ using tetrahedral elements. By using
Metis \cite{GK98}, the domain is decomposed into $N=2^k$, $k=3,4,...,9$, 
non-overlapping subdomains $\Omega_i$ with their own tetrahedral
elements; see the right plot in Fig.~\ref{fig:soldd}. 
The total number of degrees of freedom is $(2^k+1)^3$, 
$k=4,5,6,7$. We run BDDC
preconditioned GMRES iterations until the relative residual error reaches
$10^{-9}$. Three variants of BDDC preconditioners are used with corner ($C$),
corner/edge ($CE$), and corner/edge/face ($CEF$) constraints, respectively. 
The number of BDDC preconditioned GMRES iterations and the
computational time measured in seconds [s] are given in
Table~\ref{tab:the05}, with respect to the number of subdomains (row-wise)
and number of degrees of freedom (column-wise). Since the system is
unsymmetric but positive definite, the BDDC preconditioners do not show
the same typical robustness and efficiency behavior when applied to the symmetric and positive
definite system \cite{AT04}. Nevertheless, we still observe certain scalability with
respect to the number subdomains (up to $128$), in particular, with corner/edge and
corner/edge/face constraints. For $\theta=2.5$, we see improvement of BDDC
preconditioners with respect to the number of GMRES iterations and computational
time; see Table~\ref{tab:the25}. Further, we observe improved scalability with
respect to the number of subdomains as well as number of degrees of freedom.  

\begin{table}
  \begin{center}\caption{$\theta=0.5$. BDDC performance using different
      coarse level constraints ($C$/$CE$/$CEF$).}
    \begin{tabular}{ lccccccc}
      \hline
      \multicolumn{8}{c}{No preconditioner} \\ \hline
      & $8$ & $16$ & $32$ & $64$ & $128$ & $256$ & $512$ \\ \hline
      $17^3$ & $52$  & $59$  & $63$ & $71$ & $76$ & $>500$ & $>500$ \\
      &     ($0.04s$) & ($0.03s$) & ($0.03s$)& ($0.04s$) & ($0.05s$) & $-$ &$-$\\ \hline
      $33^3$ & $86$ & $100$ & $113$ & $126$ & $144$ & $>500$ & $>500$  \\ 
      & ($0.76s$) & ($0.41s$) & ($0.22s$) & ($0.18s$) & ($0.33s$) & $-$ & $-$\\ \hline
      $65^3$ &$149$ & $176$ & $206$ & $225$ & $260$ & $294$ & $>500$ \\ 
      &($35.32s$)  & ($13.55s$) & ($6.05s$) & ($2.88s$) & ($2.02s$) &($5.78s$) & $-$  \\ \hline
      $129^3$ & OoM & OoM& $346$ & $399$ & $467$ & $>500$ & $>500$ \\ 
      & ($-$) & ($-$) & ($213.23s$) & ($89.23s$) & ($41.51s$) & ($-$) & ($-$) \\ 
      \hline \\ [-0.2cm] \hline
      \multicolumn{8}{c}{$C$ (corner) preconditioner} \\ \hline
      $17^3$ & $24$ & $31$  & $30$ & $38$ & $44$ & $56$ & $114$ \\
      &  ($0.03s$)   & ($0.02s$) & ($0.02s$)& ($0.03s$) & ($0.08s$) &
      ($0.32s$) & ($2.7s$) \\ \hline
      $33^3$ & $34$ & $40$ &  $57$ & $78$ & $75$ & $247$ & $233$ \\
      & ($0.65s$) & ($0.30s$) & ($0.21s$) & ($0.19s$) & ($0.27s$) & ($2.11s$) & ($8.33s$)\\ \hline
      $65^3$ & $40$ & $64$&  $92$ & $126$ & $149$ & $166$ & $476$  \\
      & ($22.56s$) & ($10.05s$) & ($5.26s$) & ($2.95s$) & ($1.70s$) & ($4.95s$) & ($19.56s$) \\ \hline
      $129^3$ & OoM & OoM &  $173$ & $250$ & $286$ & $>500$ & $311$  \\
      & ($-$) & ($-$) & ($227.67s$) & ($113.38s$) & ($47.2s$) & ($-$) & ($20.40s$)\\ 
      \hline\\ [-0.2cm]\hline
      \multicolumn{8}{c}{$CE$ (corner+edge) preconditioner} \\ \hline
      $17^3$  & $23$ & $25$ & $24$ & $26$ & $29$ & $44$ & $70$\\ 
      &  ($0.03s$)  & ($0.02s$) & ($0.02s$) & ($0.03s$) & ($0.12s$) &
      ($0.32s$) & ($1.52s$) \\ \hline
      $33^3$ &$32$ & $32$ &  $35$ & $49$ & $40$ & $124$ & $146$\\
      & ($0.63s$)& ($0.27s$) & ($0.13s$) & ($0.15s$) & ($0.23s$) & ($1.59s$) & ($7.19s$) \\ \hline
      $65^3$ & $38$ & $56$&  $72$  & $85$ & $87$ & $89$ & $295$\\
      & ($21.93s$) & ($9.02s$) & ($4.37s$) & ($2.10s$) & ($1.23s$) & ($3.43s$) & ($21.59s$) \\ \hline
      $129^3$ & OoM & OoM &  $139$ & $177$ & $196$ & $>500$ & $218$ \\
      & ($-$ )& ($-$) & ($184.57s$) & ($82.38s$) & ($33.43s$) & ($-$) &
      ($23.76s$)\\
      \hline\\ [-0.2cm]\hline 
      \multicolumn{8}{c}{$CEF$ (corner+edge+face) preconditioner} \\ \hline
      $17^3$  &$22$ & $24$ & $23$ & $24$ & $23$ & $39$ & $72$\\ 
      &    ($0.03s$)  & ($0.02s$) & ($0.02s$) & ($0.04s$) & ($0.09$) &
      ($0.32s$) & ($1.47s$)\\ \hline
      $33^3$& $33$ & $31$ & $35$ & $46$ & $37$ & $119$ & $134$ \\ 
      & ($0.80s$)& ($0.33s$) & ($0.17s$) & ($0.17s$) & ($0.32s$) & ($2.51s$)
      & $8.2s$ \\ \hline
      $65^3$  & $39$ & $55$ & $70$ & $80$ & $81$ & $80$ & $275$\\ 
      & ($21.56s$) & ($11.52s$) & ($5.61s$) & ($2.50s$) & ($2.14s$) & ($4.04s$) & ($48.58s$)\\ \hline
      $129^3$ & OoM & OoM & $173$ & $171$ & $185$ & $>500$ & $206$ \\ 
      & ($-$) & ($-$) & ($226.93s$) & ($109.94s$) & ($45.05s$) & ($-$) & $(33.13s$) \\ 
      \hline
    \end{tabular}\label{tab:the05}
  \end{center}
\end{table}

\begin{table}\caption{ $\theta=2.5$. BDDC performance using different coarse
    level constraints ($C$/$CE$/$CEF$).}
  \begin{center}
    \begin{tabular}{ rccccccc}
      \hline
      \multicolumn{8}{c}{No preconditioner} \\ \hline
      & $8$ & $16$ & $32$ & $64$ & $128$ & $256$ & $512$ \\ \hline
      $17^3$ & $46$  & $51$  & $55$ & $61$ & $66$ & $>500$ & $>500$\\
      &($0.02s$) & ($0.02s$) & ($0.02s$)& ($0.03s$) & ($0.04s$) & $-$
      & $-$\\ \hline
      $33^3$ & $72$ & $79$ & $87$ & $99$ & $108$ & $>500$ & $>500$\\ 
      & ($0.64s$) & ($0.32s$) & ($0.15s$) & ($0.12s$) & ($0.13s$) & $-$ &  $-$ \\ \hline
      $65^3$ &$116$ & $126$ & $145$ & $163$ & $176$ & $191$ & $>500$ \\ 
      &($27.59s$)  & ($9.17s$) & ($4.07s$) & ($1.92s$) & ($1.06s$) &($2.02s$) & $-$  \\ \hline
      $129^3$ & OoM & OoM& $240$ & $271$ & $304$ & $>500$ & $382$ \\ 
      & ($-$)&($-$) & ($145.64s$) & ($58.51s$) & ($24.3s$) & ($-$) & ($12.41s$)
      \\  \hline\\ [-0.2cm]\hline
      \multicolumn{8}{c}{$C$ (corner) preconditioner} \\ \hline
      $17^3$ & $23$  & $28$  & $27$ & $32$ & $36$ & $51$ & $110$ \\
      &     ($0.02s$) & ($0.02s$) & ($0.03s$)& ($0.03s$) & ($0.07s$) &
      ($0.72s$) & ($2.48s$) \\ \hline
      $33^3$ & $30$ & $33$ & $39$ & $50$ & $50$ & $206$ & $182$ \\ 
      & ($0.60s$) & ($0.26s$) & ($0.14s$) & ($0.10s$) & ($0.09s$) & ($3.86s$) & ($6.2s$)\\ \hline
      $65^3$ &$35$ & $47$ & $61$ & $64$ & $69$ & $77$ & $287$ \\ 
      &($19.85s$)  & ($7.42s$) & ($3.47s$) & ($1.44s$) & ($0.68s$) &($1.05s$) & ($9.00s$)  \\ \hline
      $129^3$ & OoM & OoM& $94$ & $104$ & $107$ & $340$ & $112$ \\ 
      & ($-$)& ($-$) & ($124.30s$) & ($46.45s$) & ($16.90s$) & ($33.45s$) & ($5.89s$) \\ 
      \hline\\ [-0.2cm]\hline
      \multicolumn{8}{c}{$CE$ (corner+edge) preconditioner} \\ \hline
      $17^3$ & $21$  & $22$  & $22$ & $25$ & $27$ & $42$ & $80$\\
      &     ($0.03s$) & ($0.02s$) & ($0.01s$)& ($0.03s$) & ($1.78s$) &
      ($0.78s$) & ($1.65s$) \\ \hline
      $33^3$ & $27$ & $26$ & $32$ & $38$ & $32$ & $117$ & $132$ \\ 
      & ($0.54s$) & ($0.23s$) & ($0.11s$) & ($0.12s$) & ($0.15s$) & ($2.52s$) & ($5.98s$) \\ \hline
      $65^3$ &$33$ & $44$ & $51$ & $54$ & $54$ & $54$ & $235$ \\ 
      &($19.00s$)  & ($7.27s$) & ($2.98s$) & ($1.33s$) & ($0.75s$) &($1.32s$) & ($15.13s$)  \\ \hline
      $129^3$ & OoM & OoM& $82$ & $83$ & $90$ & $366$ & $94$ \\ 
      & ($-$) &($-$) & ($109.19s$) & ($38.59s$) & ($15.21s$) & ($37.00s$) & ($7.91s$) \\ 
      \hline\\ [-0.2cm]\hline
      \multicolumn{8}{c}{$CEF$ (corner+edge+face) preconditioner} \\ \hline
      $17^3$ & $21$  & $21$  & $22$ & $22$ & $22$ & $38$ & $74$ \\
      &     ($0.02s$) & ($0.01s$) & ($0.02s$)& ($0.04s$) & ($0.13s$) &
      ($0.74s$) & ($1.74s$)\\ \hline
      $33^3$ & $27$ & $26$ & $30$ & $35$ & $31$ & $115$ & $123$ \\ 
      & ($0.68s$) & ($0.28s$) & ($0.14s$) & ($0.12s$) & ($0.25s$) & ($3.79s$) & ($7.08s$) \\ \hline
      $65^3$ &$34$ & $44$ & $49$ & $52$ & $53$ & $51$ & $226$ \\ 
      &($26.64s$)  & ($9.29s$) & ($3.86s$) & ($1.61s$) & ($1.11s$) &($1.88s$) & $21.30s$  \\ \hline
      $129^3$ & OoM & OoM& $82$ & $83$ & $88$ & $369$ & $92$ \\ 
      & ($-$) & ($-$) & ($145.39s$) & ($53.10s$) & ($23.04s$) & ($52.33$) & ($13.8s$) \\ 
      \hline
    \end{tabular}\label{tab:the25}
  \end{center}
\end{table}
 
\section{Conclusions}\label{sec:cons}
In this work, we have applied two-level BDDC preconditioned GMRES methods 
%for solving the space-time finite element discretized equations of 
to the solution of finite element equations arising from the space-time discretization of a
parabolic model problem. We have compared the  performance of BDDC preconditioners with
different coarse level constraints for such an unsymmetric, but positive definite
system. The preconditioners show certain scalability 
%when we choose proper
provided that $\theta$ is sufficiently large. 
Future work will concentrate on improvement of
coarse-level corrections in order to achieve robustness with respect to
different choices of $\theta$. 

\input{referenc}
\end{document}

%% file: referenc.tex
%%%%%%%%%%%%%%%%%%%%%%%% referenc.tex %%%%%%%%%%%%%%%%%%%%%%%%%%%%%%
% sample references
% %
% Use this file as a template for your own input.
%
%%%%%%%%%%%%%%%%%%%%%%%% Springer-Verlag %%%%%%%%%%%%%%%%%%%%%%%%%%
%
% BibTeX users please use
\bibliographystyle{plain}
\bibliography{authsamp}